
\input amstex
\documentstyle{amsppt}


\topmatter
\title
On Mixing Properties of Reversible Markov Chains
\endtitle

\author
Richard C. Bradley
\endauthor

\affil
Department of Mathematics\\
Indiana University \\
Bloomington, Indiana 47405 \\
U.S.A.\\
bradleyr\@indiana.edu
\endaffil

\email
bradleyr@indiana.edu
\endemail



\abstract
It is well known that for a strictly stationary, reversible, Harris recurrent Markov chain, the $\rho$-mixing condition is equivalent to geometric ergodicity and to a ''spectral gap'' condition.  In this note, it will be shown with an example that for that class of Markov chains, the ''interlaced'' variant of the $\rho$-mixing condition fails to be equivalent to those conditions.
\endabstract

\endtopmatter

\document
\TagsOnRight
\subhead
{1.~Introduction}
\endsubhead

\medskip
This note deals with one particular aspect of the general question of where the $\rho^*$-mixing condition---the ``interlaced'' variant of the classic $\rho$-mixing condition---``fits'' into the overall scheme of mixing conditions.  Here aspects of that question will be briefly reviewed first for general strictly stationary random sequences and then for strictly stationary Markov chains.  Then we shall focus on strictly stationary Markov chains that are ``reversible'' and Harris recurrent.

For that latter class of Markov chains, it is well known that the $\rho$-mixing condition is equivalent to a ``spectral gap'' condition and also to geometric ergodicity---which in turn is equivalent to absolute regularity with exponential mixing rate.  A natural question is whether, for that particular class of Markov chains, the $\rho^*$-mixing condition can be added to that list of equivalent conditions.

In this note, we shall construct an example to show that the answer to that question is negative, even in the presence of ``nice'' extra features---a countable state space, very small marginal entropy, and information regularity with exponential mixing rate.  Aside from the various mixing conditions alluded to above, one other such condition, $\psi$-mixing, will provide useful assistance in the construction of the counterexample, even though it will not be a property of the counterexample itself.

Suppose $(\Omega, \Cal F, P)$ is a probability space.  For any two events $A$ and $B$ such that $P(A)>0$ and $P(B)>0$, define the nonnegative number
$$\eta(A,B) := {{P(A\cap B)}\over{P(A)P(B)}}\,. \tag1.1$$
For any two $\sigma$-fields $\Cal A$ and $\Cal B\subset \Cal F$, define the following four measures of dependence:  The first two are
$$\psi(\Cal A, \Cal B) := \sup_{A\in \Cal A,B\in \Cal B, P(A)>0, P(B)>0} |\eta(A,B)-1|;\,\,\hbox{and} \tag1.2$$
$$\rho(\Cal A,\Cal B) := \sup |\hbox{Corr}(f,g)| \tag1.3$$
where in (1.3) the supremum is taken over all pairs of square-integrable random variables $f$ and $g$ such that $f$ is a $\Cal A$-measurable and $g$ is $\Cal B$-measurable.  The other two measures of dependence are
$$\beta(\Cal A, \Cal B) := \sup {{1}\over{2}} \sum^I_{i=1} \sum^J_{j=1} \left| P(A_i \cap B_j) - P(A_i)P(B_j)\right|;\,\,\hbox{and} \tag1.4$$
$$I(\Cal A, \Cal B) := \sup \sum^I_{i=1} \sum^J_{j=1} \left[\eta (A_i,B_j) \log \eta (A_i,B_j)\right] P(A_i)P(B_j); \tag1.5$$
where in each of (1.4) and (1.5), the supremum is taken over all pairs of finite partitions $\{A_1,A_2,\dots, A_I\}$ and $\{B_1,B_2,\dots,B_J\}$ of $\Omega$ such that $A_i \in \Cal A$ and $P(A_i)>0$ for each $i$ and $B_j\in \Cal B$ and $P(B_j)>0$ for each $j$.  In (1.5), when necessary, $0 \log 0$ is interpreted as~$0$.  The quantity $\rho(\Cal A, \Cal B)$ in (1.3) is the ``maximal correlation'' between $\Cal A$ and $\Cal B$, and the quantity $I(\Cal A, \Cal B)$ in (1.5) is the ``coefficient of information'' between $\Cal A$ and $\Cal B$.

The factor $1/2$ in (1.4) is not important but has become standard and slightly simplifies certain inequalities.  One always has the inequalities
$$\max \{ \rho (\Cal A, \Cal B), \beta(\Cal A, \Cal B)\} \le \psi(\Cal A, \Cal B): \,\,\hbox{and} \tag1.6$$
$$[\beta(\Cal A, \Cal B)]^2\leq I(\Cal A, \Cal B) \leq (1+ \psi(\Cal A, \Cal B)) \log (1+ \psi(\Cal A, \Cal B)), \tag1.7$$
where in (1.7), when necessary, $(1+\infty) \log (1+\infty)$ is interpreted as $\infty$.  (The quantity $\psi(\Cal A, \Cal B)$, and even also $I(\Cal A,\Cal B)$, can be $\infty$.)  The second inequality in (1.7) follows trivially from (1.2) and (1.5). For the other inequalities in (1.6)--(1.7), see e.g. \cite{4, v1, Proposition 3.11(a)(b) and Theorem 5.3(I)(III)}.

Now suppose $X := (X_k, k\in \Bbb Z)$ is a strictly stationary (not necessarily Markovian) sequence of (real-valued) random variables.  For any positive integer $n$, define the following four dependence coefficients (based on (1.2)--(1.5)):
$$\split
\psi(n) &= \psi(X,n) := \psi(\sigma(X_k, k\leq 0), \sigma(X_k,k\geq n)); \\
\rho(n) &= \rho(X,n) := \rho(\sigma(X_k,k\leq 0), \sigma (X_k, k\geq n)); \\
\beta(n) &= \beta(X,n) := \beta (\sigma(X_k, k\leq 0), \sigma(X_k, k\geq n)); \,\,\hbox{and}
\endsplit $$
\vskip -.26in
$$I(n) = I(X,n) = I(\sigma(X_k,k\leq 0), \sigma(X_k,k\geq n)). \phantom{\hbox{and} } \tag1.8 $$
Here and below, the notation $\sigma(\dots)$ means the $\sigma$-field generated by $(\dots)$. For any positive integer $n$, define the dependence coefficient $$\rho^*(n) = \rho^*(X,n) := \sup \rho(\sigma(X_k,k\in S), \sigma(X_k, k\in T)) \tag1.9$$
where the supremum is taken over all pairs of nonempty disjoint sets $S,T\subset \Bbb Z$ such that $\hbox{dist}(S,T):= \min_{s\in S,t\in T} |s-t| \ge n$.  (The two index sets $S$ and $T$ can be ``interlaced'', with each one having elements between ones in the other set.)  Obviously $\psi(1) \geq \psi(2) \geq \psi(3) \geq\dots$, and the analogous comment applies to each of the other dependence coefficients in (1.8) and (1.9). The strictly stationary sequence $X$ is said to satisfy ``$\psi$-mixing'' if $\psi(n) \to 0$ as $n\to \infty$, ``$\rho$-mixing'' if $\rho(n) \to 0$ as $n\to \infty$, ``absolute regularity'' (or ``$\beta$-mixing'') if $\beta(n) \to 0$ as $n\to \infty$, ``information regularity'' if $I(n) \to 0$ as $n\to \infty$, and ``$\rho^*$-mixing'' if $\rho^*(n) \to 0$ as $n\to \infty$.

The $\rho^*$-mixing condition can be regarded as an ``interlaced'' variant of $\rho$-mixing.  For strictly stationary sequences, corresponding ``interlaced'' variants of $\psi$-mixing, absolute regularity, and information regularity are equivalent to $m$-dependence and are not treated further here; see e.g. \cite{4, v1, Theorem 5.11 and Remark 5.12(II)}.

The information regularity condition is of interest in information theory; see e.g. the book by Pinsker \cite{22}. The other mixing conditions formulated above have each played a substantial role in the development of limit theory for weakly dependent random sequences.  Among the $\psi$-mixing, $\rho$-mixing, absolute regularity, and information regularity conditions---the conditions here that involve just ``past'' and ``future''---the $\psi$-mixing condition implies the other three by (1.6)--(1.7), information regularity implies absolute regularity by (1.7), and (see e.g. \cite{4, v1, Chart 5.22}) there are no other implications.  However, it is not known whether or not $\psi$-mixing implies $\rho^*$-mixing.  That is part of the ongoing question of precisely where the $\rho^*$-mixing condition ``fits'' into the broader scheme of mixing conditions.

Here is some motivation for that question. The strictly stationary sequences that satisfy $\rho^*$-mixing (formulated with different terminology), with exponential mixing rate, were one of the contexts in which Charles Stein \cite{26} introduced what later became known as ``Stein's method'' for obtaining bounds on rates of convergence in central limit theorems and other related theorems for weakly dependent sequences. (For more on Stein's method, see e.g. the book by Chen, Goldstein, and Shao \cite{5}.)  Later, Goldie and Greenwood \cite{8} and others used the $\rho^*$-mixing condition in more conventional ways in central limit theory for random sequences and random fields.  The author \cite{1} showed that for some central limit theory under just finite second moments, the $\rho^*$-mixing condition works effectively even with an arbitrary slow mixing rate.  That theme was developed further in various papers such as Miller \cite{15} and Tone \cite{28}. The same theme was devel
 oped by Peligrad \cite{20}, and later also in other papers such as Utev and Peligrad \cite{29} and Tone \cite{27}, in a modified form with $\rho^*$-mixing replaced by the weaker pair of conditions (i)~$\rho^*(n) <1$ for some $n\geq 1$ and (ii)~(Rosenblatt) ``strong mixing'' (or ``$\alpha$-mixing'') with again an arbitrarily slow mixing rate.  Because of the results in some of those papers, if an affirmative answer were established for the (still open) question of whether (for strictly stationary sequences) the classic ``$\phi$-mixing'' condition implies $\rho^*$-mixing, then as an immediate corollary, one would have an affirmative answer to two (still unsolved) long-standing conjectures in central limit theory for $\phi$-mixing sequences: Ibragimov's Conjecture (see Ibragimov and Linnik \cite{9, p.~393, Problem (3)}) and a related conjecture of Iosifescu \cite{10, p.~56} (the one involving the functional central limit theorem). (The best results on those two conjectures so f
 ar, giving an affirmative answer to at least the spirit of both conjectures, are those of Peligrad \cite{18}\cite{19}.) All that motivates the ongoing question of precisely where $\rho^*$-mixing fits into the scheme of mixing conditions.

Now let us review that question in the narrower context of strictly stationary Markov chains. It is well known that for a given strictly stationary Markov chain $X := (X_k, k\in \Bbb Z)$ and a given positive integer $n$, the dependence coefficients in (1.8) satisfy the following equalities:
$$\split
\psi(n) &= \psi(X,n) = \psi(\sigma(X_0), \sigma(X_n));\\
\rho(n) &= \rho(X,n) = \rho(\sigma(X_0), \sigma(X_n));\\
\beta(n) &= \beta(X,n) = \beta(\sigma(X_0), \sigma(X_n));\,\,\hbox{and}
\endsplit
$$
\vskip -.25in
$$I(n) = I(X,n) = I(\sigma(X_0),\sigma(X_n)). \phantom{\,\, \hbox{and}}\tag1.10
$$
(See e.g.\ \cite{25} or \cite{4, v1, Theorem 7.3}.) For a given strictly stationary Markov chain $X := (X_k$, $k\in \Bbb Z)$ and a given pair of positive integers $m$ and $n$, one has the following well known inequalities:
$$
\split
\psi(m+n) &\leq \psi (m) \cdot \psi(n); \\
\rho(m+n) &\leq \rho(m) \cdot \rho(n); \,\,\, \hbox{and}
\endsplit$$
\vskip-.26in
$$\rho^*(m+n) \leq \rho^*(m) \cdot \rho^*(n). \phantom{: \,\,\, \hbox{and} } \tag1.11$$
For the first two inequalities, apply e.g.\ (1.10) and \cite{4, v1, Theorem 7.4(a)(c)} (to the ``Markov triplet'' of $\sigma$-fields $(\sigma(X_0), \sigma(X_m), \sigma(X_{m+n})))$; and for a similar proof of the third inequality, adapt e.g.\ the proof of \cite{4, v1, Theorem 7.5(I)(a$'$)}. Of course by (1.11) and induction, for strictly stationary Markov chains, the mixing rates for $\psi$-mixing, $\rho$-mixing, and $\rho^*$-mixing must be exponential. (By ``exponential mixing rate'', we mean of course, for example for $\psi$-mixing, only that there exists $r\in (0,1)$ such that $\psi(n) = O(r^n)$ as $n\to \infty$. That allows $\psi(n)$ to be ``large'', even $\infty$, for ``small'' $n$, and it allows $\psi(n) \to 0$ ``more quickly than exponentially fast'', even $m$-dependence.) For strictly stationary Markov chains, for the absolute regularity and information regularity conditions, the mixing rate need not be exponential. See e.g.\  the examples in Davydov \cite{7} or \cite{
 4, v1, Example 7.11}. By a result of Nummelin and Tuominen \cite{16, Theorem 2.1}, building on earlier work of Nummelin and Tweedie \cite{17}, a strictly stationary Markov chain satisfies the classic ``geometric ergodicity'' condition if and only if it satisfies absolute regularity with exponential mixing rate.  For more on geometric ergodicity, see e.g.\ Meyn and Tweedie \cite{14} or \cite{4, v2, Theorem 21.19}. Now Rosenblatt \cite{25, p.~214, line~$-3$ to p.~215, line~16} constructed a family of strictly stationary Markov chains which satisfy $\rho$-mixing but not absolute regularity. A slightly modified version of those examples, again a strictly stationary Markov chain, was explicitly shown in \cite{4, v1, Example 7.16} to satisfy $\rho^*$-mixing but not absolute regularity. The author \cite{2, Theorem 1.2} showed that for strictly stationary Markov chains, $\psi$-mixing implies $\rho^*$-mixing. (The result there was more general in that (i)~it allowed non-stationarity,
  with appropriate modifications of the formulations of dependence coefficients such as in (1.8) and (1.10), and (ii)~it assumed only the ``lower half'' of the $\psi$-mixing condition.) As an immediate corollary (also treated separately in \cite{4, v1, Theorem 7.15} with a more gentle proof), every strictly stationary, finite-state, irreducible, aperiodic Markov chain is $\rho^*$-mixing.

The author \cite{3, p.~719, Theorem~1 and p.~725, Remark~2} constructed a class of strictly stationary, countable-state Markov chains that satisfy $\rho$-mixing, absolute regularity, and information regularity with mixing rates that can be made arbitrarily fast (short of $m$-dependence) such that the $\rho^*$-mixing condition fails to hold.  Trivially those Markov chains are irreducible and aperiodic, but they do not satisfy the condition of ``reversibility'', to which we now turn.

A given strictly stationary (real-valued) Markov chain $X := (X_k, k\in \Bbb Z)$ is said to be ``reversible'' if the distribution (on $\Bbb R^{\Bbb Z})$ of the ``time-reversed'' Markov chain $(X_{-k}, k\in \Bbb Z)$ is the same as that of $X$ itself---equivalently if the distribution (on $\Bbb R^2)$ of the random vector $(X_1,X_0)$ is the same as that of $(X_0,X_1)$.  The strictly stationary Markov chain in \cite{4, v1, Example 7.16} alluded to above---which satisfies $\rho^*$-mixing but not absolute regularity---is easily seen to be reversible, but it is not Harris recurrent.

Now let us focus on the class of strictly stationary (real-valued) Markov chains that are both reversible and Harris recurrent.  For that particular class of Markov chains, it is well known, from works such as Roberts and Rosenthal \cite{23}, Roberts and Tweedie \cite{24, p.~39}, and Kontoyannis and Meyn \cite{11} \cite{12, Proposition~1.2}, that the $\rho$-mixing condition is equivalent to geometric ergodicity (absolute regularity with exponential mixing rate), equivalent to a certain ``spectral gap'' condition, and (consequently---see e.g.\ Longla and Peligrad \cite{13, Theorem 4}) equivalent to the condition $\rho(1)<1$. A natural question is whether, for that particular class of Markov chains, the $\rho^*$-mixing condition can be added to that list of equivalent conditions.

The answer turns out to be negative, even in the presence of very nice extra conditions involving a countable state space, small marginal entropy, and information regularity. That will be shown in Theorem 1.1 below, after a little more notation is given.

For a given purely atomic $\sigma$-field $\Cal A$ with (finitely many or countably many) atoms $A_1, A_2, A_3,\dots$, the ``entropy of $\Cal A$'' is defined by
$$H(\Cal A) := - \sum_k \left[ P(A_k) \cdot \log P(A_k)\right] = I(\Cal A, \Cal A). \tag1.12$$
(The first equality is the definition; the second is a trivial by-product of (1.5)---see e.g.\ \cite{22} or \cite{4, v1, Section~5.9}.) For any purely atomic $\sigma$-field $\Cal A$ and any (not necessarily atomic) $\sigma$-field $\Cal B$, one has the elementary inequality
$$I(\Cal A, \Cal B) \leq H(\Cal A). \tag1.13$$
(Again see \cite{22} or \cite{4, v1, Section~5.9}.)

Here is the main result of this note:

\proclaim{Theorem 1.1} Suppose $r$ is a number such that $0<r<1$. Then there exists a strictly stationary, countable-state, irreducible, aperiodic Markov chain $X:= (X_k, k\in \Bbb Z)$ which is reversible and has the following three properties:
\roster
\item One has that $H(\sigma(X_0)) \leq r$.
\item For every positive integer $n$,
$$\max\left\{ \rho(X,n), \beta(X,n), I(X,n)\right\} \le r^n. \tag1.14$$
\item For every positive integer $n$,
$$\rho^*(X,n) = \rho(\sigma(X_0), \sigma(X_{-n},X_n)) = 1. \tag1.15$$
\endroster
\endproclaim

By (1.7) there is obviously some ``redundancy'' in (1.14). With the state space being countable, the properties ``irreducible'' and ``aperiodic'', and also ``recurrent'', are of course trivial consequences of (1.14).  In (1.15), the failure of $\rho^*$-mixing is manifested on pairs of index sets ($\{0\}$ and $\{-n,n\}$) of smallest possible cardinality.

Theorem 1.1 will be proved in Section 3.  Section 2 will provide a special class of strictly stationary, finite-state, irreducible, aperiodic Markov chains that will be used as ``building blocks'' for the construction in Section~3.

\subhead 2. The Building Blocks
\endsubhead

\medskip
In this section, a special class of Markov chains will be developed.  They will be used in Section~3 in the construction of the Markov chain $X$ for Theorem~1.1.

The following inequality will be needed:
$$\forall\,\, N\in \{3,4,5,\dots\,\}, \,\,\,\forall\,\, \varepsilon \in (0, 1/3], \quad 1-2 \cdot \left[ \varepsilon^{2N-1} + \sum^{N-1}_{u=1} \varepsilon^{2u}\right] > {{1}\over{2}}. \tag2.1$$

\definition{Definition 2.1}
Suppose $N\geq 3$ is an integer and $\varepsilon \in (0,1/3]$.  A strictly stationary Markov chain $Y:= (Y_k, k\in \Bbb Z)$ is said to satisfy ``Condition ${\Cal S}(N, \varepsilon)$'' if the following three conditions hold:
\roster
\item The state space of $Y$ is the set $\{0,1,2,\dots,N\}$.
\item The (marginal) distribution $\mu_{N,\varepsilon, i} := P(Y_0 =i)$, $ i\in \{0,1,\dots,N\}$ of the random variable $Y_0$ is as follows (see (2.1)):
$$\split
&\mu_{N,\varepsilon,0} = \varepsilon^2 + 1-2 \cdot \left[ \varepsilon^{2N-1} + \sum^{N-1}_{u=1} \varepsilon^{2u}\right]; \\
&\forall\, m\in \{1,2,\dots,N-2\},\,\, \mu_{N,\varepsilon,m} = \varepsilon^{2m} +\varepsilon^{2(m+1)};\\
&\mu_{N,\varepsilon,N-1}= \varepsilon^{2(N-1)} + \varepsilon^{2N-1}; \,\,\,\hbox{and}\\
&\mu_{N,\varepsilon,N}= \varepsilon^{2N-1}.
\endsplit
\tag2.2$$

\item The distribution of the random vector $(Y_0,Y_1)$ is as follows:
$$\split
&P((Y_0,Y_1) =(0,0)) = 1-2 \cdot \left[ \varepsilon^{2N-1} + \sum^{N-1}_{u=1} \varepsilon^{2u}\right]; \\
&\forall\,\, m \in \{1,2,\dots,N-1\},\\
&\qquad P((Y_0,Y_1) = (m-1,m)) = P((Y_0,Y_1) = (m,m-1)) = \varepsilon^{2m};\\
&P((Y_0,Y_1) =(N-1,N)) = P((Y_0,Y_1) = (N,N-1)) = \varepsilon^{2N-1}; \,\,\hbox{and} \\
&P((Y_0,Y_1) = (i,j)) = 0 \,\,\,\hbox{if $i=j\geq 1$ or $|i-j| \geq 2$.}
\endsplit
\tag2.3
$$
\endroster
\enddefinition

(In various notations such as $\mu_{N,\varepsilon,i}$ in (2.2), it will be handy to retain explicitly the parameters $N$ and $\varepsilon$ for clarity.)

It is easy to check (see (2.1)) that (i)~in each of (2.2) and (2.3), the assigned probabilities are nonnegative and add up to 1, (ii)~eq.\ (2.3) is compatible with (2.2), and (iii)~from (2.3) the distribution of the random vector $(Y_0,Y_1)$ is the same as that of $(Y_1,Y_0)$, and hence the Markov chain $Y$ is reversible.

For a given integer $N\ge 3$ and a given $\varepsilon \in (0,1/3]$, for the Markov chain $Y := (Y_k, k\in \Bbb Z)$ in Definition~2.1 satisfying Condition $\Cal S(N,\varepsilon)$, the one-step transition probabilities will be denoted by $p_{N,\varepsilon,i,j} := P(Y_1 = j \mid Y_0 =i)$ for $i,j\in \{0,1,\dots, N\}$, and for any positive integer $m$, the $m$-step transition probabilities will be denoted by $p^{(m)}_{N,\varepsilon,i,j} := P(Y_m = j\mid Y_0 =i)$.

In what follows, the notation $a_\varepsilon \sim b_\varepsilon$ as $\varepsilon \to 0+$, for positive numbers $a_\varepsilon$ and $b_\varepsilon$, means that $\lim_{\varepsilon\to 0+} a_\varepsilon/b_\varepsilon =1$. By (2.2), for each integer $N\geq 3$, one has the following:
$$\split
&\mu_{N,\varepsilon, 0} \rightarrow 1 \,\,\,\hbox{as }\,\, \varepsilon \to 0+; \,  \\
&\forall\,\, m\in \{1,2,\dots,N-1\},\,\,\mu_{N,\varepsilon,m} \sim \varepsilon^{2m}\,\,\hbox{as $\varepsilon \to 0+$; and}\\
&\mu_{N,\varepsilon,N} \sim \varepsilon^{2N-1}\,\,\,\hbox{as}\,\,\,\varepsilon \to 0+.
\endsplit
\tag2.4
$$
Also, by (2.2) and (2.3) and trivial arithmetic, for each integer $N\geq 3$, one has the following:
$$\split
&p_{N,\varepsilon,0,0} \rightarrow 1 \,\,\,\hbox{as}\,\,\,\varepsilon \to 0+;\\
&\forall\, m\in \{1,2,\dots,N\}, \,\, p_{N,\varepsilon, m,m-1} \to 1 \,\,\hbox{as}\,\,\varepsilon \to 0+; \\
&\forall\, m\in \{0,1,\dots,N-2\}, \, p_{N,\varepsilon, m,m+1} \sim \varepsilon^2 \,\,\hbox{as}\,\,\varepsilon \to 0+;\\
&p_{N,\varepsilon, N-1,N} \sim \varepsilon \,\,\hbox{as}\,\,\varepsilon \to 0+;\,\,\hbox{and}\\
&\forall\, \varepsilon \in (0,1/3],\, p_{N,\varepsilon,i,j} =0 \,\,\hbox{if}\,\,i=j\geq 1\,\,\hbox{or}\,\,|i-j|\geq 2.
\endsplit
\tag2.5
$$
Finally, by (2.5) and a simple argument, for each integer $N\geq 3$, one has the following:
$$\split
&\forall\, m\geq 1,\, p^{(m)}_{N,\varepsilon,0,0} \rightarrow 1\,\,\hbox{as}\,\,\varepsilon \to 0+;\\
&\forall\, m\in\{1,2,\dots,N-1\}, \, p^{(m)}_{N,\varepsilon,0,m} \sim \varepsilon^{2m} \,\,\hbox{as}\,\,\varepsilon \to 0+;\\
&p^{(N)}_{N,\varepsilon,0,N} \sim \varepsilon^{2N-1} \,\,\hbox{as}\,\, \varepsilon \to 0+;\\
&\hbox{if}\,\,N\geq j>i\geq 0\,\,\hbox{then}\,\,p^{(j-i)}_{N,\varepsilon,j,i} \rightarrow 1\,\,\hbox{as}\,\, \varepsilon \to 0+;\,\,\hbox{and}\\
&\hbox{if} \,\, 1\leq j\leq N\,\, \hbox{ and}\,\, m\geq j\,\, \hbox{ then}\,\,p^{(m)}_{N,\varepsilon,j,0} \rightarrow 1\,\,\hbox{as}\,\,\varepsilon \to 0+.\endsplit
\tag2.6
$$

The following lemma lists the properties of the Markov chains in Definition 2.1 that will be employed when those Markov chains are used in Section~3 to construct the Markov chain $X$ for Theorem~1.1.

\proclaim{Lemma 2.2} Suppose $N\geq 3$ is an integer. For each $\varepsilon \in (0,1/3]$, let $Y^{\varepsilon)} = Y^{(N,\varepsilon)} := (Y^{(\varepsilon)}_k$, $k\in \Bbb Z)$ be a strictly stationary Markov chain that satisfies Condition $\Cal S(N,\varepsilon)$.  The the following seven statements hold:
\roster
\item For each $\varepsilon \in (0,1/3]$, the Markov chain $Y^{(\varepsilon)}$ is reversible.
\item One has that $\lim_{\varepsilon \to 0+} P(Y^{(\varepsilon)}_0 =0) =1$.
\item One has that $\lim_{\varepsilon\to 0+} H(\sigma(Y^{(\varepsilon)}_0)) =0$.
\item One has that $\lim_{\varepsilon \to 0+} \rho(Y^{(\varepsilon)},1) =0$.
\item For every integer $m$ such that $1\leq m<N/2$, one has that
$$\lim_{\varepsilon \to 0+} \rho\left(\sigma(Y^{(\varepsilon)}_0), \sigma\left( Y^{(\varepsilon)}_{-m}, Y^{(\varepsilon)}_m\right)\right) = 1 \tag2.7$$
\item One has that $\lim_{\varepsilon \to 0+} \psi(Y^{(\varepsilon)},5N) =0$.
\item For any $r\in (0,1)$, there exists a number $\varepsilon_{N,r} \in (0,1/3]$ such that
$$\forall\,\,\varepsilon \in (0,\varepsilon_{N,r}],\,\, \forall  \,\,n\in \Bbb N,\,\,\, I(Y^{(\varepsilon)},n) \leq r^n. \tag2.8$$
\endroster
\endproclaim

{\bf Proof.} Throughout this proof, the integer $N\geq 3$ is fixed.  Properties (1) and (2) (in Lemma 2.2) hold by (2.3) and (2.4) respectively (as was already noted).

{\bf Proof of Property (3).}  For each $\varepsilon \in(0,1/3]$, the atoms of the $\sigma$-field $\sigma(Y^{(\varepsilon)}_0)$ are the events $\{Y^{(\varepsilon)}_0 =j\}$, $j\in \{0,1,\dots,N\}$; and hence by (1.12) (see (2.2)),
$$H\left( \sigma(Y^{(\varepsilon)}_0)\right) = - \sum^N_{j=0} \mu_{N,\varepsilon,j} \log \mu_{N,\varepsilon,j}.$$
By (2.4), $\mu_{N,\varepsilon,0}\to 1$ as $\varepsilon \to0+$ and for every $j\in \{1,2,\dots,N\}$, $\mu_{N,\varepsilon,j} \to 0$ as $\varepsilon\to 0+$. Since $x\log x \to 0$ as either $x\to 0+$ or $x\to 1$, Property (3) holds.

{\bf Proof of Property (4).} To simplify the argument, we shall consider pairs of events instead of pairs of square-integrable random variables.  For any two $\sigma$-fields $\Cal A$ and $\Cal B \subset \Cal F$, define the measure of dependence
$$\lambda(\Cal A, \Cal B) := \sup_{A\in \Cal A, B\in \Cal B} {{|P(A\cap B)-P(A)P(B)|}\over{[P(A)P(B)]^{1/2}}}\,\,. \tag2.9$$
Here and in such fractions below, $0/0$ is interpreted as $0$.  For any $\alpha \in (0,1)$, there exists $\delta \in (0,1)$ such that whenever $\Cal A$ and $\Cal B$ are $\sigma$-fields satisfying $\lambda (\Cal A, \Cal B) \leq \delta$, one has that $\rho(\Cal A,\Cal B) \leq \alpha$. A reasonably elementary proof of this fact can be found in \cite{4, v1, Theorem 4.15}. For a very sharp version of this fact, with a more complicated proof, see Peyre \cite{21}. Here we just note that to prove Property (4), it suffices to prove that
$$\lambda\left(\sigma(Y^{(\varepsilon)}_0), \sigma(Y^{(\varepsilon)}_1)\right) \longrightarrow 0\,\,\hbox{as}\,\,\varepsilon \to 0+. \tag2.10$$

For any $\varepsilon \in (0,1/3]$ and any two sets $\Gamma,\Lambda \subset \{0,1,\dots,N\}$, define the quantity
$$f(\varepsilon,\Gamma,\Lambda) := {{|P(\{Y^{(\varepsilon)}_0 \in \Gamma\}\cap\{Y^{(\varepsilon)}_1\in \Lambda\}) - P(Y^{(\varepsilon)}_0 \in \Gamma)\cdot P(Y^{(\varepsilon)}_1\in \Lambda)|}\over{[P(Y^{(\varepsilon)}_0 \in \Gamma)\cdot P(Y^{(\varepsilon)}_1\in \Lambda)]^{1/2}}}\,\,, \tag2.11$$
(where again $0/0$ is interpreted as $0$). Of course for any $\varepsilon \in (0,1/3]$ and any $k\in \Bbb Z$, the events in the $\sigma$-field $\sigma(Y^{(\varepsilon)}_k)$ are precisely the ones $\{Y^{(\varepsilon)}_k \in \Gamma\}$ for $\Gamma \subset \{0,1,\dots,N\}$. Hence to prove (2.10) it suffices to prove that
$$\left[\sup_{\Gamma,\Lambda \subset \{0,1,\dots,N\}} f(\varepsilon,\Gamma,\Lambda)\right] \longrightarrow 0 \,\,\hbox{as}\,\,\varepsilon \to 0+. \tag2.12$$

Now by the trivial equality $[P(F^c \cap G) - P(F^c)P(G)] = -[P(F\cap G) - P(F) P(G)]$ for events $F$ and $G$, eq.\ (2.9) does not change if the supremum there is taken over just the events $A\in \Cal A$ and $B\in \Cal B$ such that $P(A) \leq 1/2$ and $P(B) \leq 1/2$. For a given $\varepsilon \in (0,1/3]$ and a given $k\in \Bbb Z$, one has that $P(Y^{(\varepsilon)}_k =0) > 1/2$ by (2.1) and (2.2), and hence for a given set $\Gamma \subset \{0,1,\dots,N\}$, $P(Y^{(\varepsilon)}_k \in \Gamma) \leq 1/2$ if and only if $0 \notin \Gamma$. Hence (see (2.11)) the ``goal equation'' (2.12) does not change if the supremum there is taken over just the sets $\Gamma, \Lambda \subset \{1,2,\dots,N\}$ (i.e.\ with the state $0$ excluded).

Also, there are only finitely many subsets of $\{1,2,\dots,N\}$.  Consequently, one can set up the proof of (2.12) as follows: Let $\Gamma$ and $\Lambda$ be arbitrary fixed sets $\subset \{1,2,\dots,N\}$. To prove (2.12) (and thereby Property (4)), it suffices to show for this $\Gamma$ and $\Lambda$ that
$$\lim_{\varepsilon \to 0+} f(\varepsilon,\Gamma,\Lambda) =0. \tag2.13$$

If either $\Gamma$ and $\Lambda$ is empty, then $f(\varepsilon, \Gamma,\Lambda) =0$ for all $\varepsilon \in (0,1/3]$ and (2.13) holds trivially.  Therefore, assume that neither $\Gamma$ and $\Lambda$ is empty.

Define the integers $g,\ell \in \{1,2,\dots,N\}$ by $g:= \min \Gamma$ and $\ell := \min \Lambda$. The proof of (2.13) will be divided into the three cases $g<\ell$, $g>\ell$, and $g =\ell$.  By symmetry (e.g.\ in (2.11)) and reversibility, the arguments for the first two cases are exactly analogous.  It will suffice to consider just the two cases $g<\ell$ and $g=\ell$.

{\it Case 1:} $g<\ell$. Then $1\leq g\leq N-1$ and $2\leq \ell\leq N$. By (2.4) and a simple calculation, as $\varepsilon \to 0+$ one has that $P(Y^{(\varepsilon)}_0 \in \Gamma) \sim \varepsilon^{2g}$ and $P(Y^{(\varepsilon)}_1 \in \Lambda) \sim \varepsilon^{2\ell}$ (resp. $\varepsilon^{2N-1}$) if $2\leq \ell \leq N-1$ (resp. $\ell=N$). Hence in any case
$$\left[ P(Y^{(\varepsilon)}_1 \in \Lambda)/P(Y^{(\varepsilon)}_0 \in \Gamma)\right]^{1/2} \longrightarrow 0 \,\,\hbox{as}\,\,\varepsilon \to 0+. \tag2.14$$
Also, for any two events $F$ and $G$ such that $0<P(F) \leq P(G)$, the fractions $P(F\cap G)/[P(F)P(G)]^{1/2}$ and $P(F)P(G)/[P(F)P(G)]^{1/2}$ are each trivially bounded above by $[P(G)/P(F)]^{1/2}$. Hence for a given $\varepsilon \in (0,1/3]$ the right hand side of (2.11) is bounded above by the left hand side of (2.14). Hence by (2.14), eq.\ (2.13) holds (in the case $g<\ell$).

 {\it Case 2:} $g=\ell$.  Consider first the subcase where $g=\ell=N$. Then $\Gamma=\Lambda=\{N\}$. By (2.3), for any $\varepsilon \in(0,1/3]$, $P(Y^{(\varepsilon)}_0 = Y^{(\varepsilon)}_1 = N)=0$ and hence by (2.11), $f(\varepsilon,\Gamma, \Lambda)$ equals $[P(Y^{(\varepsilon)}_0 =N) \cdot P(Y^{(\varepsilon)}_1 =N]^{1/2}$.  That converges to $0$ as $\varepsilon \to 0+$ by (2.4). Thus (2.13) holds if $g=\ell=N$.

Now finally consider the case where $g=\ell \in \{1,2,\dots,N-1\}$. By (2.3), $P(Y^{\varepsilon)}_0 = Y^{(\varepsilon)}_1 =g) =0$ and hence
$$P\left(\{Y^{(\varepsilon)}_0 \in \Gamma\} \cap \{Y^{(\varepsilon)}_1 \in \Lambda\}\right) \leq P\left(Y^{(\varepsilon)}_0 \geq g+1\right) + P\left( Y^{(\varepsilon)}_1 \geq g+1\right). \tag2.15$$
By (2.4), the right hand side of (2.15) is $\sim 2\varepsilon^{2g+2}$ (resp. $2\varepsilon^{2N-1}$) if $g+1 \leq N-1$ (resp. $g+1=N$). Also by (2.4) and a simple calculation,
$$\left[ P(Y^{(\varepsilon)}_0 \in \Gamma)\cdot P(Y^{(\varepsilon)}_1 \in \Lambda)\right]^{1/2} \sim \varepsilon^{2g}\,\,\hbox{as}\,\,\varepsilon \to 0+. \tag2.16$$
Hence (whether $g+1\leq N-1$ or $g+1=N$)
$$[\hbox{LHS of (2.15)}] = o([\hbox{LHS of (2.16)}]) \,\,\hbox{as}\,\,\varepsilon \to 0+.$$
Also of course by (2.16) itself,
$$[\hbox{LHS of (2.16)}]^2 = o([\hbox{LHS of (2.16)}]) \,\,\hbox{as}\,\,\varepsilon \to 0+.$$
Hence by (2.11), eq.\ (2.13) holds (if $g=\ell \in \{1,\dots,N-1\})$. That completes the proof of Property (4).

{\bf Proof of Property (5).} Let $m$ be an arbitrary fixed positive integer such that $m< N/2$. For each $\varepsilon \in (0,1/3]$, define the events
$$A_\varepsilon := \{Y^{(\varepsilon)}_0 =N\} \,\,\,\hbox{and}\,\,\,B_\varepsilon := \{Y^{(\varepsilon)}_{-m} = Y^{(\varepsilon)}_m = N-m\}. \tag2.17$$
To complete the proof of Property (5), it suffices to show that
$$\lim_{\varepsilon \to 0+} \hbox{Corr}(I(A_\varepsilon), I(B_\varepsilon)) = 1 \tag2.18$$
(Here and below, $I(\dots)$ denotes the indicator function.)

By (2.17) and (2.4),
$$P(A_\varepsilon) \sim \varepsilon^{2N-1}\,\,\hbox{as}\,\, \varepsilon \to 0+. \tag2.19$$
Also, for each $\varepsilon \in (0,1/3]$, by (2.17), the Markov property, and reversibility (Property (1)), $P(B_\varepsilon \mid A_\varepsilon) = [p^{(m)}_{N,\varepsilon,N,N-m}]^2$. That converges to 1 as $\varepsilon \to 0+$ by (2.6). Hence by (2.19),
$$P(A_\varepsilon \cap B_{\varepsilon}) \sim \varepsilon^{2N-1} \,\,\hbox{as}\,\,\varepsilon \to 0+. \tag2.20$$

To complete the proof of (2.18) (and of Property (5)), it will suffice to show that
$$P(B_\varepsilon) \sim \varepsilon^{2N-1}\,\,\hbox{as}\,\,\varepsilon \to 0+. \tag2.21$$
For then by (2.19), (2.20), and (2.21), one will have that as $\varepsilon \to 0+$, $P(A_\varepsilon \cap B_\varepsilon )\sim [P(A_\varepsilon) \cdot P(B_\varepsilon)]^{1/2}$ and $P(A_\varepsilon) \cdot P(B_\varepsilon) = o([P(A_\varepsilon) \cdot P(B_\varepsilon)]^{1/2})$ and hence Cov$(I(A_\varepsilon), I(B_\varepsilon)) \sim [P(A_\varepsilon) \cdot P(B_\varepsilon)]^{1/2}$, and by (2.19) and (2.21) and a trivial calculation, one will also have that Var~$I(A_\varepsilon) \sim P(A_\varepsilon)$ and Var~$I(B_\varepsilon) \sim P(B_\varepsilon)$ as $\varepsilon \to 0+$, and (2.18) will follow.

By (2.20), to prove (2.21), it suffices to show that
$$P(A^c_\varepsilon \cap B_\varepsilon) = o(\varepsilon^{2N-1}) \,\,\,\hbox{as}\,\,\varepsilon \to 0+. \tag2.22$$

Let $S$ denote the set of all $(2m+1)$-tuples
$${\bold j} := (j_{-m}, j_{-m+1},\dots, j_m) \in \{0,1,\dots,N\}^{2N+1} \tag2.23$$
such that
$$j_{-m} = j_m = N-m \,\,\,\hbox{and}\,\,\,j_0 \not= N. \tag2.24$$
For each $\varepsilon \in (0,1/3]$ and each ${\bold j}\in S$, using the representation of $\bold j$ in (2.23), define the event
$$C_{\varepsilon,{\bold j}} := \bigcap^m_{k=-m} \left\{Y^{(\varepsilon)}_k =  j_k\right\}. \tag2.25$$
Then by (2.17), for each $\varepsilon \in (0,1/3]$, $P(A^c_\varepsilon \cap B_\varepsilon) = \sum_{{\bold j} \in S} P(C_{\varepsilon, {\bold j}})$.  Since the set $S$ is finite, the rest of the proof of (2.22) can be set up as follows:

Let ${\bold j} \in S$ be arbitrary but fixed.  To prove (2.22) (and thereby Property (5)), it suffices to prove for this $\bold j$ that
$$P(C_{\varepsilon,{\bold j}}) = o(\varepsilon^{2N-1}) \,\,\,\hbox{as}\,\,\, \varepsilon \to 0+. \tag2.26$$

Refer again to the representation of $\bold j$ in (2.23). If $| j_{k+1} -  j_k| \geq 2$ for some $k\in \{-m,-m+1,\dots,m-1\}$, then $P(C_{\varepsilon,j}) =0$ for each $\varepsilon \in (0,1/3]$ by (2.25) and (2.5) and we are done.  Therefore, assume that for each $k\in \{-m,-m+1,\dots,m-1\}$, $|j_{k+1} -  j_k| \leq 1$.

Because of that last assumption, one now has that
$$\forall\,\,k\in \{-m,-m+1,\dots,m\},\,\,\, 1\leq  j_k \leq N-1. \tag2.27$$
Here the second inequality $ j_k \leq N-1$ holds by (both parts of) (2.24), and the first inequality $1\leq { j}_k$ holds by the assumption $m< N/2$ (which implies $N-m>m$) and the double equality in (2.24).

If $ j_k =  j_{k+1} \geq 1$ for any $k \in \{-m,-m+1,\dots,m-1\}$, then $P(C_{\varepsilon,{\bold j}}) =0$ for each $\varepsilon \in (0,1/3]$ by (2.25) and (2.5) and we are done.  Therefore (see (2.27)) we now assume that $| j_{k+1} -  j_k| =1$ for every $k\in \{-m, -m+1,\dots m-1\}$.  Combining that with (2.27) and the double equality in (2.24), one now has the following:

The difference $ j_{k+1} - j_k$ equals $+1$ for exactly $m$ values of $k\in \{-m,-m+1,\dots,m-1\}$, and equals $-1$ for the other $m$ such values. Hence by (2.25), (2.27), (2.4), and (2.5),
$$\split
P(C_{\varepsilon,{\bold j}}) &= \mu_{N,\varepsilon,N-m} \cdot \prod^{m-1}_{k=-m} p_{N,\varepsilon,  j(k), j(k+1)} \\
&\sim \varepsilon^{2(N-m)} \cdot (\varepsilon^2)^m \cdot 1^m \\
&= \varepsilon^{2N}\,\,\,\hbox{as}\,\,\varepsilon \to 0+
\endsplit$$
(where $ j(k)$ and $ j(k+1)$ mean $ j_k$ and $ j_{k+1}$).  Hence (2.26) holds.  That completes the proof of Property (5).

{\bf Proof of Property (6).}  Refer to the notations in (1.1), (2.2), and (2.6). For any $\varepsilon \in (0,1/3]$ and any pair of elements $i,j \in \{0,1,\dots,N\}$, define the number
$$\split
g(\varepsilon,i,j) &:= \eta\left(\{Y^{(\varepsilon)}_0 =i\}, \{Y^{(\varepsilon)}_{5N} = j\}\right) \\
&= {{ P\left( Y^{(\varepsilon)}_{5N} = j\mid Y^{(\varepsilon)}_0 =i\right)}\over{P\left(Y^{(\varepsilon)}_{5N} = j\right)}} = {{p^{(5N)}_{N,\varepsilon,i,j}}\over{\mu_{N,\varepsilon,j}}} \,.
\endsplit
\tag2.28$$
Of course for any $\varepsilon \in (0,1/3]$ and any integer $k$, the atoms of the $\sigma$-field $\sigma(Y^{(\varepsilon)}_k)$ are precisely the events $\{Y^{(\varepsilon)}_k =i\}$ for $i\in \{0,1,\dots,N\}$.  By an elementary argument (see e.g.\ \cite{4, v1, Proposition 3.21(4)}), for any two $\sigma$-fields $\Cal A$ and $\Cal B$ that are purely atomic, eq.\ (1.2) remains unchanged if the supremum there is taken only over the atoms $A$ and $B$ of $\Cal A$ and $\Cal B$ respectively.  Hence by (1.10), to prove that $\psi(Y^{(\varepsilon)}, 5N) \to 0$ as $\varepsilon\to 0+$, it suffices to prove that $g(\varepsilon,i,j) \to 1$ uniformly in $i,j \in \{0,1,\dots,N\}$ as $\varepsilon \to 0+$. Since the state space $\{0,1,\dots,N\}$ is finite, that task can be set up as follows:

Let $i$ and $j$ each be an arbitrary fixed element of $\{0,1,\dots,N\}$. To prove Property (6), it suffices to prove for this $i$ and $j$ that
$$g(\varepsilon,i,j) \rightarrow 1 \,\,\hbox{as}\,\,\varepsilon \to 0+. \tag2.29$$

Consider first the case $i=j=0$. By (2.4) and (2.6), $\mu_{N,\varepsilon,0} \to 1$ and $p^{(5N)}_{N,\varepsilon,0,0} \to 1$ as $\varepsilon \to 0+$.  Hence by (2.28), eq.\ (2.29) holds.

Next consider the case where $i \in \{1,2,\dots,N\}$ and $j=0$. Again by (2.4) and (2.6), $\mu_{N,\varepsilon,0} \to 1$ and $p^{(5N)}_{N,\varepsilon,i,0} \to 1$ as $\varepsilon \to 0+$; and hence by (2.28), eq.\ (2.29) holds.

Now for the case where $i=0$ and $j\in \{1,2,\dots,N\}$, one has by reversibility (Property (1)) and the first equality in (2.28) that $g(\varepsilon,0,j) = g(\varepsilon,j,0)$ for each $\varepsilon \in(0,1/3]$.  Thus in this case, (2.29) holds as a corollary of the preceding case.

Now finally suppose that $i,j\in \{1,2,\dots,N\}$.  To complete the proof of Property (6), what remains is to prove (2.29) for this case.  From (2.4) and the last term in (2.28), one sees that in order to prove (2.29), our task is to verify that
$$\hbox{as}\,\, \varepsilon \to 0+, \quad p^{(5N)}_{N,\varepsilon,i,j} \sim \cases \varepsilon^{2j} &\text{if $1\leq j\leq N-1$} \\
\varepsilon^{2N-1} &\text{if $j=N$.} \endcases
\tag2.30$$

Let $H$ denote the set of all $(5N+1)$-tuples
$${\bold h} := (h_0,h_1,h_2,\dots,h_{5N}) \in \{0,1,\dots,N\}^{5N+1} \tag2.31$$
such that
$$h_0 =i\,\,\,\hbox{and}\,\,\,h_{5N}=j. \tag2.32$$
For each $\varepsilon \in (0,1/3]$ and each ${\bold h} \in H$, using the representation (2.31)--(3.23), define the quantity
$${\bold p}_\varepsilon ({\bold h}) := \prod^{5N}_{u=1} p_{N,\varepsilon,h(u-1),h(u)} \tag2.33$$
(where $h(u-1)$ and $h(u)$ mean $h_{u-1}$ and $h_u$). Then for each $\varepsilon \in (0,1/3]$,
$$p^{(5N)}_{N,\varepsilon,i,j} = \sum_{{\bold h}\in H} {\bold p}_\varepsilon ({\bold h}). \tag2.34$$

Now using the representation (2.31)--(2.32) for elements $\bold h$ of $H$, we shall define three subsets $H_1,H_2,H_3 \subset H$ such that $H\supset H_1 \supset H_2 \supset H_3$.

First, let $H_1$ denote the set of all ${\bold h} \in H$ such that
$$\forall\,\,\, k\in \{1,2,\dots,5N\},\,\,\hbox{ either}\,\,|h_k-h_{k-1}| =1 \,\,\,\hbox{or}\,\,\, h_{k-1} = h_k =0. \tag2.35$$

Next, let $H_2$ denote the set of all ${\bold h} \in H_1$ (i.e.\ satisfying (2.35)) such that there exists $k\in \{1,2,\dots, 5N-1\}$ for which $h_k =0$. If $\bold h$ is an element of $H_2$ and $\ell$ is the greatest element of $\{1,2,\dots, 5N-1\}$ such that $h_\ell =0$, then by (2.35) (and (2.32)) the ordered pairs $(h_{u-1},h_u)$, $u \in \{\ell+1,\ell+2,\dots,5N\}$ must include $(0,1),(1,2),\dots,(j-1,j)$, each at least once.

Finally, let $H_3$ denote the set of all ${\bold h}\in H_2$ such that among all of the ordered pairs $(h_{u-1}, h_u)$, $u\in \{1,2,\dots, 5N\}$, there are exactly $j$ of the form $(q-1,q)$, $q\in \{1,2,\dots,N\}$. If $\bold h$ is an element of $H_3$ and $\ell$ is the greatest element of $\{1,2,\dots,5N-1\}$ such that $h_\ell =0$, then the following hold: (i)~By (2.32), (2.35), and a simple argument, $\ell = 5N-j$ must hold, and the ordered pairs $(h_{u-1}, h_u)$, $u = \ell+1, \ell+2,\dots, 5N$ are respectively $(0,1),(1,2),\dots,(j-1,j)$. (That exhausts the quota of the $j$ ordered pairs of the form $(q-1,q)$.)  (ii)~Hence by (2.32) and (2.35), the ordered pairs $(h_{u-1},h_u)$, $u=1,2,\dots,i$ are respectively $(i,i-1),(i-1,i-2),\dots,(1,0)$; and the ordered pairs $(h_{u-1},h_u)$, $i+1 \leq u \leq 5N-j$ are all $(0,0)$. It now follows that the only element of $H_3$ is the $(5N+1)$-tuple
$${\bold h}^* := (i,i-1,\dots,1,0,0,\dots,0,1,2,\dots,j). \tag2.36$$
Now by (2.5),
$$p_{N,\varepsilon,0,1} \cdot p_{N,\varepsilon,1,2} \cdot \dots \cdot p_{N,\varepsilon,j-1,j} \sim [\text{RHS of (2.30)] as $\varepsilon \to 0+$.} \tag2.37$$
Hence by (2.33), (2.36), and (2.5), ${\bold p}_\varepsilon ({\bold h}^*) \sim$ [RHS of (2.30)] as $\varepsilon \to 0+$.

Consequently by (2.34) (and the fact that $H$ is a finite set), to prove (2.30) it suffices to show that for all ${\bold h} \in H-H_3$ (i.e.\ all ${\bold h}\ \in H$ except ${\bold h}^*$),
$${\bold p}_\varepsilon ({\bold h}) = o([\text{RHS of (2.30)}]) \,\,\hbox{as}\,\,\varepsilon \to 0+. \tag2.38$$

First suppose ${\bold h} \in H_2-H_3$. Then (see the definition of $H_3$), among the ordered pairs $(h_{k-1},h_k)$, $1\leq k\leq 5N$, there are at least $j+1$ of the form $(q-1,q)$, $q\in \{1,2,\dots,N\}$, and (see the second sentence after (2.35)), those ordered pairs include $(0,1),(1,2),\dots,(j-1,j)$ each at least once. The product in (2.33) includes all of the terms in the left side of (2.37) as well as at least one ``extra'' term of the form $p_{N,\varepsilon,q-1,q}$, $q \in \{1,2,\dots,N\}$.  By (2.5), that ``extra'' term is either $\sim \varepsilon^2$ or $\sim \varepsilon$ as $\varepsilon \to 0+$.  Combining that with (2.37) (and the fact that all terms in the product in (2.33) are bounded above by 1), one has by (2.33) the desired result (2.38).

Next suppose ${\bold h} \in H_1-H_2$. Then by (2.35) and the definition of $H_2$, one has that $|h_k -h_{k-1}| =1$ for every $k\in \{1,2,\dots,5N\}$. There must be at least $2N$ indices $k\in \{1,2,\dots,5N\}$ such that $h_k-h_{k-1} =1$ (for otherwise---see (2.32)---in the telescoping sum $j-i = \sum^{5N}_{k=1} (h_k -h_{k-1})$, there would be fewer than $2N$ terms $+1$, hence more than $3N$ terms $-1$, hence $j-i < -N$, contradicting the stipulation $i,j \in \{1,2,\dots,N\})$.  Thus by (2.5) there would be at least $2N$ indices $k$ for which the term $p_{N,\varepsilon,h(k-1),h(k)}$ in the product (2.33) is either $\sim \varepsilon^2$ or $\sim \varepsilon$ as $\varepsilon \to 0+$. Hence by (2.33), one again has the desired result (2.38).

Finally suppose ${\bold h} \in H-H_1$. Then (see (2.35) in the definition of $H_1$) there exists $k\in \{1,2,\dots, 5N\}$ such that $|h_k-h_{k-1}| \geq 2$ or $h_{k-1} = h_k \geq 1$. For that $k$, $p_{N,\varepsilon, h(k-1),h(k)} =0$ for each $\varepsilon \in (0,1/3]$ by (2.5). Hence by (2.33), again eq.\ (2.38) holds.

Thus (2.38) holds for all ${\bold h} \in H$ other than ${\bold h}^*$. That completes the proof of Property (6).

{\bf Proof of Property (7).} Suppose $0<r<1$.  Applying Properties (3) and (6), let $\varepsilon_{N,r} \in (0,1/3]$ be such that
$$\forall\, \varepsilon(0,\varepsilon_{N,r}],\,\, H(\sigma(Y^{(\varepsilon)}_0)) \leq r^{5N}\,\,\,\hbox{and}\,\,\,\psi(Y^{(\varepsilon)},5N) \le (r/2)^{10N}. \tag2.39$$

Now let $\varepsilon \in (0,\varepsilon_{N,r}]$ be arbitrary but fixed.  Our task is to prove for this $\varepsilon$ that $I(Y^{(\varepsilon)},n)\leq r^n$ for all $n\in \Bbb N$.

In the case where $1\leq n\leq 5N$, one has by (1.10), (1.13), and (2.39) that
$$I(Y^{(\varepsilon)},n) =I(\sigma(Y^{(\varepsilon)}_0), \sigma(Y^{(\varepsilon)}_n))
\leq H(\sigma(Y^{(\varepsilon)}_0)) \leq r^{5N} <r^n,
$$
giving the desired result.

Now suppose instead that $n>5N$.  Let $M$ be the positive integer such that $M\cdot (5N) < n\leq (M+1)\cdot(5N)$.  By (1.11) and induction, followed by (2.39),
$$\psi(Y^{(\varepsilon)},M\cdot 5N) \leq [\psi(Y^{(\varepsilon)},5N)]^M \leq (r/2)^{10MN}.$$
Hence by (1.7) and the inequality $\log (1+x) \leq x$ for $x>0$,
$$\split
I(Y^{(\varepsilon)},n) &\leq I(Y^{(\varepsilon)},M\cdot 5N) \\
&\leq (1+\psi(Y^{(\varepsilon)},M\cdot 5N)) \cdot \log (1+\psi(Y^{(\varepsilon)}, M\cdot 5N)) \\
&\leq (1+ (r/2)^{10MN}) \cdot \log (1+(r/2)^{10MN}) \\
&\leq 2\cdot (r/2)^{10MN} \\
&\leq r^{2M\cdot 5N} \leq r^{(M+1)\cdot 5N} \leq r^n,
\endsplit$$
again giving the desired result.  That completes the proof of Property (7), and of Lemma 2.2.

\subhead
{3. Proof of Theorem 1.1}
\endsubhead

\medskip
First a key lemma is stated. Let $\Bbb N$ denote the set of all positive integers.

\proclaim{Lemma 3.1} Suppose $\Cal A_n$ and $\Cal B_n$, $n\in \Bbb N$ are $\sigma$-fields $\subset \Cal F$, and the $\sigma$-fields $\Cal A_n \vee \Cal B_n$, $n\in \Bbb N$ are independent. Then
$$\rho\left( \bigvee^\infty_{n=1} \Cal A_n, \bigvee^\infty_{n=1} \Cal B_n\right) = \sup_{n\in \Bbb N} \rho(\Cal A_n, \Cal B_n); \,\,\hbox{and} \tag3.1$$
$$I\left( \bigvee^\infty_{n=1} \Cal A_n, \bigvee^\infty_{n=1} \Cal B_n \right) = \sum_{n\in \Bbb N} I(\Cal A_n, \Cal B_n). \tag3.2$$
\endproclaim

Eq.\ (3.1) is due to Cs\'aki and Fischer \cite{6, Theorem 6.2}. Eq.\ (3.2) is a classic fact from information theory (see e.g.\ \cite{22}). Proofs of both (3.1) and (3.2) can be found e.g.\ in \cite{4, v1, Theorems 6.1--6.2}.

\medskip
{\bf Proof of Theorem 1.1.} As in the statement of that theorem, suppose $0<r<1$.
Applying Lemma 2.2, for each integer $N\geq 3$, let $Z^{(N)} := (Z^{(N)}_k, \, k\in \Bbb Z)$ be a strictly stationary Markov chain that is reversible, has state space $\{0,1,\dots,N\}$, and satisfies the following five conditions:
$$P(Z^{(N)}_0 =0) \geq 1-2^{-N}; \tag3.3$$
$$H(\sigma(Z^{(N)}_0)) \leq 2^{-N}r; \tag3.4$$
$$\rho(Z^{(N)},1) \leq r; \tag3.5$$
$$I(Z^{(N)},n) \leq (2^{-N}r)^{2n}\,\,\hbox{for every $n\in \Bbb N$; and} \tag3.6$$
$$\rho\left( \sigma(Z^{(N)}_0), \sigma(Z^{(N)}_{-m},Z^{(N)}_m)\right) \geq 1 - {{1}\over{N}}\,\,\hbox{for every $m\in \Bbb N$ such that $m<N/2$.} \tag3.7$$
Let these sequences $Z^{(3)}, Z^{(4)}, Z^{(5)},\dots$ be constructed in such a way that they are independent of each other.

By strict stationarity, (3.3), and the Borel--Cantelli Lemma, for each $k\in \Bbb Z$,
$$P\left(Z^{(N)}_k \not= 0 \,\,\hbox{for infinitely many} \,\, N\geq 3\right) =0.$$
Changing the definition of the sequences $Z^{(N)}$, $N\geq 3$ on a set of probability $0$ if necessary, we henceforth assume that for every $\omega \in \Omega$ and every integer $k$, $Z^{(N)}_k(\omega) \not= 0$ for at most finitely many $N\geq 3$.

Let $S$ denote the set of all sequences $h:= (h_3,h_4,h_5,\dots)$ of nonnegative integers such that (i)~$h_N \in \{0,1,\dots,N\}$ for each $N\geq 3$ and (ii)~$h_N \not= 0$ for at most finitely many integers $N\geq 3$. This set $S$ is countable.

Define the sequence $X:= (X_k$, $k\in \Bbb Z)$ of $S$-valued random variables as follows: For each $k\in \Bbb Z$,
$$X_k := \left( Z^{(3)}_k, Z^{(4)}_k, Z^{(5)}_k,\dots\right). \tag3.8$$

By standard arguments, this sequence $X$ is a strictly stationary Markov chain.  Also, since each of the Markov chains $Z^{(N)}$, $N\geq 3$ is reversible, it follows by a simple argument that the Markov chain $X$ is reversible.

Property (1) in Theorem~1.1 holds since by (3.8) and (3.3),
$$H(\sigma(X_0)) = \sum^\infty_{N=3} H\left( \sigma(Z^{(N)}_0)\right) \leq r.$$
The equality here is a standard property of entropy (and can be checked here trivially from (1.12) and (3.2)).

Now let us verify Property (2) in Theorem~1.1.  By (3.8), (3.1), and (3.5), $\rho(X,1) = \sup_{N\geq 3} \rho(Z^{(N)},1) \leq r$.  Hence by (1.11) and induction,
$$ \forall\, n\in \Bbb N,\quad \rho(X,n) \leq r^n. \tag3.9$$
Next, for each $n\in \Bbb N$, by (3.8), (3.2), and (3.6),
$$
I(X,n) = \sum^\infty_{N=3} I(Z^{(N)},n) \leq \sum^\infty_{N=3} (2^{-N}r)^{2n}
< \sum^\infty_{N=3} 2^{-N} r^{2n} <r^{2n}.
\tag3.10$$
Hence by (1.7), $\beta(X,n)< r^n$ for every $n\in \Bbb N$.  From that and (3.9) and (3.10), Property (2) in Theorem~1.1 has been verified.

Now let us verify Property (3).  Suppose $n\in \Bbb N$.  For each integer $N\geq 2n+1$, by (3.8) and (3.7),
$$\rho(\sigma(X_0), \sigma(X_{-n}, X_n)) \geq \rho(\sigma(Z^{(N)}_0), \sigma(Z^{(N)}_{-n}, Z^{(N)}_n)) \geq 1- 1/N.$$
Property (3) follows.

It was already noted that the state space $S$ of the Markov chain $X$ is countable.  As a trivial consequence of that and Property (2), the Markov chain $X$ is irreducible and aperiodic.  All features of $X$ asserted in Theorem~1.1 have been verified; the proof is complete.

\remark{Remark 3.2} Given the interest here in the absolute regularity condition (because of its connection to geometric ergodicity), one might want to obtain the inequality $\beta(X,n) \leq r^n$, $n\in \Bbb N$ in Property (2) in Theorem~1.1 by an argument more direct than through information regularity and (1.7). To do that, one would first, in the proof of Lemma~2.2, slightly modify the proof of its Property (7) to obtain $\beta(Y^{(\varepsilon)},n) \leq r^n$ in place of $I(Y^{(\varepsilon)},n) \leq r^n$ in (2.8). In that argument for Property (7), one would replace the first inequality in (2.39) by $P(Y^{(\varepsilon)}_0 =0) \geq 1 -r^{5N}/2$, and then obtain $\beta(Y^{(\varepsilon)},n) \leq r^{5N}$ (say for $1\leq n\leq 5N)$ from a well known property of $\beta(\Cal A, \Cal B)$ for $\sigma$-fields $\Cal A$ and $\Cal B$ when (say) $\Cal A$ has a ``big atom'' (see e.g.\ \cite{4, v1, Proposition 3.19}). Also in that argument for Property (7), one would employ (1.6) instead o
 f the second inequality in (1.7). Then in the proof of Theorem~1.1, eq.\ (3.6) would be replaced by $\beta(Z^{(N)},n) \leq (2^{-N}r)^n$, and one would then derive the inequality $\beta(X,n) \leq r^n$ from a well known inequality for $\beta(\dots)$ similar to eq.\ (3.2) in the context of Lemma~3.1 (see e.g.\ \cite{4, v1, Theorem 6.2}).
\endremark

\NoBlackBoxes


\end